\documentclass[12pt,english]{article}
\usepackage{amsfonts,here,graphicx}
\usepackage{cite,amsfonts,amssymb} 
\def\qed{\kern8pt \vrule height5pt depth0pt width5pt}

\def\cF{{\cal F}}

\def\cH{{\cal H}}
\def\cK{{\cal K}}

\def\cM{{\cal M}}

\def\cP{{\cal P}}

\def\sup{\mathop{\rm sup}}

\def\proof{\noindent \medskip {\bf Proof:}$\;\;$}

\def\ass#1#2\endass{\vskip5pt plus2pt \noindent{\bf (A.#1)} {\it #2} \vskip5pt
plus2pt }

\newtheorem{prop}{Proposition}

\newtheorem{lem}{Lemma}

\newtheorem{theor}{Theorem}

\newtheorem{cor}{Corollary}

\title{\bf Function--valued stochastic convolutions arising in 
integrodifferential equations}

\author{\large\sf Anna Karczewska \\
 \\
 Institute of Mathematics,
 University of Zielona G\'ora\\
 ul. Podg\'orna 50, 65-246 Zielona G\'ora, Poland\\
 e-mail: A.Karczewska@im.uz.zgora.pl\\
}
\date{\today}

\begin{document}

\maketitle

\def\thefootnote{}

\begin{abstract}
We study stochastic convolutions providing by fundamental solutions of 
a class of integrodifferential 
equations which interpolate the heat and the wave equations. 
We give sufficient
condition for the existence of function--valued convolutions in terms of the
covariance kernel of a noise given by spatially homogeneous Wiener process.
\end{abstract}

\section{Introduction}\label{int}

The paper is concerned with the stochastic integrals of the form
\begin{equation} \label{eq1}
 \int_0^t P(t-s) * (b(u(s))dW(s)), \quad t\in [0,T]\;.
\end{equation}
In (\ref{eq1}), $P$ is a fundamental solution of some integrodifferential 
equation introduced in the next section, 
 $u$ is a stochastic process specified next,
$W$ is a spatially homogeneous Wiener process and 
$b$: $\mathbb{R}\rightarrow\mathbb{R}$ denotes the random field.

The integral (\ref{eq1}) looks out formally like those considered in \cite{PeZa2}
where the existence of function--valued solutions of nonlinear stochastic wave 
and heat equations have been analyzed. This is the well-known fact (see, e.g.\
\cite{Mi}) that solutions to nonlinear wave and heat equations may be represented
in terms of fundamental solutions to the Cauchy problems for  
equations, wave or heat, respectively. 
This idea has been used in \cite{PeZa2} for representing
solutions of stochastic nonlinear equations. In both obtained formulas, the
stochastic integrals of the form (\ref{eq1}) have appeared. 
Unfortunately, no analogous formula (with fundamental solution) exists for
the solution to nonlinear Volterra equation.

As we have already mentioned, the paper  \cite{PeZa2} is concerned with 
function--valued solutions to the stochastic nonlinear wave and heat equations
and provides necessary and sufficient conditions for the existence of such
solutions. The problem of existence of function--valued solutions to linear and
nonlinear stochastic wave and heat equations has been investigated by many people.
They obtained several conditions in terms of function coefficients, the covariance
kernel or spectral measure of the noise $W$. We refer to some papers only:  
\cite{Da}, \cite{DaFr}, \cite{DaMu}, \cite{KaZa1}, \cite{KaZa2}, \cite{Pe}, 
\cite{PeZa1} and \cite{PeZa2}, 
for more information see references therein.
But only a few papers are concerned with the stochastic Volterra equation. The
paper \cite{KaZa3} deals with linear Volterra equation and is written in the
spirit of the above mentioned papers.

The aim of this paper is to provide conditions under which the stochastic
integral  (\ref{eq1}) is well--defined process with values in the space 
$L_v^2 = L^2(\mathbb{R},v(x)dx)$, where $v$ is some test
function. 
We follow the idea used in some part of \cite{PeZa2} for our study on the
interval $[-R,R]$. It could be done because the fundamental 
solutions $P_\alpha$, $1<\alpha <2$,
considered in the paper are similar in some sense to fundamental solution of the
wave equation. Next part of our studies, particularly the passing to the limit
as $R\rightarrow +\infty$, bases mostly on specific properties of solutions 
$P_\alpha$.
The paper proposes a framework for a study of nonlinear 
stochastic Volterra equations.
Till now, to the best of our knowledge, there has been appeared no paper
concerning  function--valued solutions to nonlinear Volterra equations. We hope
the results obtained in the paper will be the first step towards this direction.

In the paper we study only the one-dimensional case. To treat regularity of the
stochastic integral (\ref{eq1}) in the cases $d>1$ by this method 
will require development
of the results analogous to Fujita's ones for the higher dimensions. 
And, to date, this work has not been done.

\section{Deterministic integrodifferential equations}

The following integrodifferential equation
\begin{equation} \label{eq2}
 u(t,x) = g(x) + \int_0^t h(t-s)\Delta u(s,x)\,ds \;,
\end{equation}
$t>0$, $x\in \mathbb{R}$, has been treated by many authors (see, e.g.\
\cite{FrSh}, \cite{Fr},  \cite{Fu},\cite{Pr} and references therein). 
The equation (\ref{eq2}) is usually used to describe the heat conduction in
materials with memory. Very important results have been obtained (see \cite{Fu}
and \cite{ScWy}) in the case 
$ h(t) = {t^{\alpha-1}}/{\Gamma(\alpha)}$, for $1\leq\alpha\leq 2$,
where $\Gamma$ is the gamma function. With such a kernel the equation (\ref{eq2})
reads
\begin{equation} \label{eq3}
 u(t,x) = g(x) + \frac{1}{\Gamma(\alpha)} \int_0^t (t-s)^{\alpha-1}
 \Delta u(s,x)\,ds, \quad  1\leq\alpha\leq 2 \;.
\end{equation}
We can see that the family of equations (\ref{eq3}) with 
$1\leq\alpha\leq 2$ interpolates, in some sense, 
the heat equation (when $\alpha =1$) and the wave equation 
(when $\alpha =2$). Both papers \cite{Fu} and \cite{ScWy} provide the
representation of the solutions to the equations (\ref{eq3}) and they are
complementary to one another. 
Additionally, Fujita has characterized the solutions to (\ref{eq3}) in detail. 
Because we shall use his results in the paper, we first 
recall some facts from Fujita's work.

Let $S(\mathbb{R})$ denote the space of rapidly decreasing functions and 
$C([0,+\infty );S(\mathbb{R}))$ be the space consisting of $S(\mathbb{R})$--valued
continuous functions on $[0,+\infty )$. For $1\leq\alpha\leq 2$ and $t>0$, we
define the function
$$
 q_\alpha (t,\xi) := 
 \exp [-t|\xi|^\delta \,e^{-i\pi\gamma \mbox{\tiny sgn} (\xi)/2} ]\;,
\mbox{~~where~~} \delta = 2/\alpha \mbox{~~and~~} \gamma = 2-2/\alpha.
$$

Define $P_\alpha(t,x)$ as follows
\begin{equation} \label{eq4}
 P_\alpha(t,x) :=  \frac{1}{2\pi} \int_{-\infty}^{+\infty} q_\alpha (t,\xi)
 e^{-ix\xi} d\xi\;.
\end{equation}

The function $P_\alpha(t,x)$ has the following properties
\begin{equation} \label{eq5}
\left\{\begin{array}{ll}
  P_\alpha(t,x) \geq 0\;, & ~~t\in (0,+\infty), ~~ x\in \mathbb{R}  \\
  ~ & ~ \\
  \int\limits_{-\infty}^{+\infty} P_\alpha(t,x) dx = 1\;, & ~~t\in (0,+\infty)
\end{array} \right. 
\end{equation}

i.e.\ is a probability density function and
\begin{equation} \label{eq6}
  P_\alpha(t,x) = P_\alpha(xt^{-\alpha/2})\, t^{-\alpha/2},\quad t\in (0,+\infty),
  ~~ x\in \mathbb{R}\;, 
\end{equation}
where $P_\alpha(x)= P_\alpha(1,x)$.\\[2mm]

Let us recall the representation of the solution to the equation (\ref{eq3}).
We assume that the function $g$ in (\ref{eq3}) belongs to the space 
$S(\mathbb{R})$.
\begin{theor} \label{th.2.1} (Theorem A, \cite{Fu})\\
 For every $1\leq\alpha\leq 2$, the equation  (\ref{eq3}) has a unique solution 
 $u_\alpha (t,x)$ given by
$$  u_\alpha (t,x) =\left\{
\begin{array}{ll}
 \frac{1}{\alpha} \int\limits_{-\infty}^{+\infty} g(x-y) P_\alpha(t,|y|)\, dy,
   & ~~\mbox{for~~} 1\leq\alpha <2 \\
   ~ & ~ \\
 \frac{1}{2}\,[g(x+t)+g(x-t)]\,,  & ~~\mbox{for~~} \alpha =2\;.
\end{array}
\right.       $$
\end{theor}

Hence, $\frac{1}{\alpha}\,P_\alpha(t,|x|) ~\mbox{with}~ 1\leq\alpha\leq 2$ 
is the fundamental solution of the
equation (\ref{eq3}). This means that this function is the integral kernel of
the operator acting from the initial data $g$ to the solution $u_\alpha (t,x)$
of (\ref{eq3}).

The next theorems provide important properties of the fundamental solution to
the equation (\ref{eq3}).

Let us define, for $1\leq \alpha<2$ and $x\in \mathbb{R}$
 the following functions
\begin{eqnarray} \label{eq7}
 a_\alpha(x) &=& |x|^{\frac{2}{\alpha}}\exp \left[ \frac{\pi i}{\alpha} 
 \mbox{sgn} (x)\right] \;,\nonumber\\
  &&\nonumber\\
 b_\alpha(x) &=& |x|^{\frac{2}{\alpha}}\exp \left[ -\frac{\pi i}{\alpha} 
 \mbox{sgn} (x)\right] \;,
\end{eqnarray}
and
$$  f_\alpha(x) = \left \{
  \begin{array}{ll} \frac{\sin (\alpha\pi)}{\pi} \int\limits_0^\infty
  \frac{x^2\,t^{\alpha-1}\,e^{-t}}{t^{2\alpha}+2x^2 t^\alpha 
  \cos (\alpha\pi)+x^4} \quad\quad & (x\neq 0) \\
  1-\frac{2}{\alpha}\quad\quad  & (x= 0) 
  \end{array} \right. \;.
$$
\vspace{3mm}
\begin{theor} \label{th.2.2a} (Lemma 1.4, \cite{Fu})\\
For every $1\leq \alpha<2$
\begin{equation} \label{eq8}
 \cF^{-1}[F_\alpha](x) = \frac{1}{\alpha}\,P_\alpha(|x|)\;,
\end{equation}
where $$
F_\alpha(x) = \frac{1}{\alpha}\{\exp (a_\alpha(x)) + \exp (b_\alpha(x))\}
  + f_\alpha(x) \quad\quad \mathrm{for} \quad x\in \mathbb{R}  $$
and functions $a_\alpha, b_\alpha, f_\alpha$ are defined like in (\ref{eq7}).
\end{theor}
\vspace{3mm}
\begin{cor} \label{cor1} 
 From (\ref{eq8}) we obtain
$$  \cF^{-1}[\exp (a_\alpha)](x) = P_\alpha(-x) \quad \mathrm{and} \quad
 \cF^{-1}[\exp (b_\alpha)](x) = P_\alpha(x)\;.
$$ 
\end{cor}

\newpage
\begin{theor} \label{th.2.2} (Theorem B, \cite{Fu})\\
 For $1<\alpha < 2$ we have: 
 \begin{enumerate}
  \item $P_\alpha(t,|x|)$ is continuous for $t\in (0,+\infty)$, 
    $x\in \mathbb{R}$.
  \item $P_\alpha(t,|x|)$ takes its extreme values at 
    $x=\pm c_\alpha\,t^{\alpha/2}$ (maximum) and $x=0$ (minimum), where
    $c_\alpha >0$ is a constant determined by $\alpha$.
    The solution is monotone elsewhere.
  \item  $P_\alpha(t,|x|)$ never vanishes for $t\in (0,+\infty)$, 
    $x\in \mathbb{R}$.
 \end{enumerate}
\end{theor}

\vspace{3mm}
\noindent{\bf Comment:} {\em 
By Theorem \ref{th.2.2} we see, that the fundamental solution to  
(\ref{eq3}) has the similar property to that of the wave equation. For 
both equations, the points, where the fundamental solution takes its maximum,
propagate with finite speed.}

\vspace{2mm}
The below picture illustrates Theorem \ref{th.2.2}.
\begin{figure}[h]
\label{fig1}
\begin{center}
 \resizebox{0.95\textwidth}{!}{\includegraphics{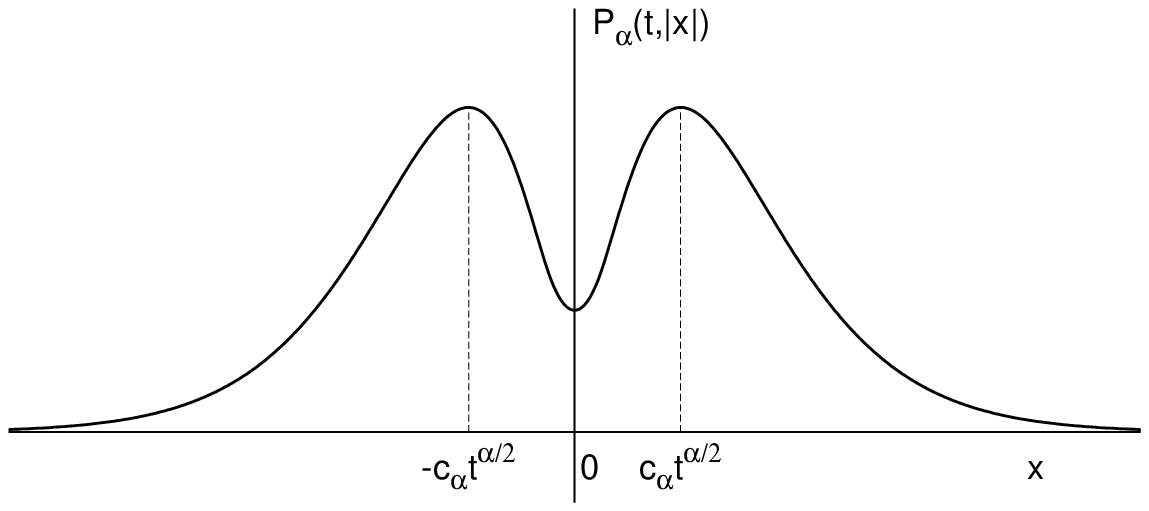}}
\end{center}
\end{figure}

The fundamental solution to the equation (\ref{eq3}) is well-known for 
the limiting case $\alpha=0$, not considered here,
$\alpha=1$, that is for the heat equation  and $\alpha=2$, 
that is for the wave equation.
For some more information we refer to \cite{Fu}, \cite{ScWy} and \cite{Pr}.
Let us notice that the point of view on the equation (\ref{eq3}) represented
by the authors is different. While Friedman and Pr\"uss emphasize the
representation of the fundamental solution $P_\alpha(t,|x|)$ to (\ref{eq3})
for general case $0<\alpha<2$ in terms of the Mittag-Leffler functions 
(Friedman \cite{Fr} was the first who observed that fact), Fujita and 
Schneider with Wyss join the fundamental solution with $\alpha$-stable (or
L\'evy) probability distributions. Moreover, the latter authors propose a
possible physical application of fractional diffusion ($0<\alpha\leq 1$)
not considered in this paper.

As we have already written, we consider the same family of the equations 
(\ref{eq3}) like Fujita, that is for  $1\leq \alpha\leq 2$.

\section{Stochastic integral}\label{stoi}

Because several papers (cf.\ \cite{PeZa1}, \cite{KaZa2}, \cite{KaZa3},
\cite{PeZa2}) contain detailed description of the stochastic integral used 
in this paper, we recall only the most important facts indispensable for
understanding the remaining part of the paper.

Let $S(\mathbb{R}^d)$ denote the space of all so called rapidly decreasing
functions on $\mathbb{R}^d$ and $S'(\mathbb{R}^d)$ be the space of tempered
distributions on $\mathbb{R}^d$ (see, e.g.\ \cite{Yo}).
Denote by $S_{(s)}(\mathbb{R}^d)$ the space of functions 
$\varphi\in S(\mathbb{R}^d)$ such that $\varphi=\varphi_{(s)}$, where 
$\varphi_{(s)}(x)=\overline{\varphi(-x)},~~x\in \mathbb{R}^d$. 
Analogously, we
denote by $S'_{(s)}(\mathbb{R}^d)$ the space of all distributions 
$\xi\in S'(\mathbb{R}^d)$ such that $\langle \xi,\varphi\rangle =
\langle \xi,\varphi_{(s)}\rangle$ for $\varphi \in S(\mathbb{R}^d)$.
(In the whole paper, the value of a distribution $\xi$ on a test function 
$\varphi$ will be denoted by $\langle \xi,\varphi\rangle$.)

Assume that $(\Omega,\cF,(\cF_t)_{t\geq 0},\mathbb{P})$ 
is a complete probability space. As we
have already written, the noise process $W$ is a spatially homogeneous Wiener
process. Such process has been used by many authors, see e.g.\
the papers mentioned above 
and references therein.
Shortly speaking, $W$ is an $S'(\mathbb{R}^d)$-valued Wiener process having the
following properties: 
\begin{enumerate}
 \item for every $\varphi \in S(\mathbb{R}^d)$, 
  $\langle W(t),\varphi\rangle_{t\in [0,+\infty)}$ is a real-valued process;
  \item process $W$ has the covariance of the form 
  $E\langle W(t),\varphi\rangle\langle W(t),\psi\rangle =
  s \wedge t\,\langle \Gamma, \varphi * \psi_{(s)}\rangle$,
  where $\varphi \in S(\mathbb{R}^d)$, $\psi \in S_{(s)}(\mathbb{R}^d)$
  and $\Gamma$ is a positive-definite distribution in $S'(\mathbb{R}^d)$.
\end{enumerate}

Let us emphasize that the spatially homogeneous Wiener process $W$ for any fixed
$t\geq 0$ becomes a stationary, Gaussian, generalized random field.

Since  $\Gamma$  is a positive-definite tempered distribution, there exists a
positive symmetric so called {\em tempered measure}
$\mu$ on $\mathbb{R}^d$ such that 
$\Gamma=\cF (\mu)$.
(Let us recall that the measure  $\mu$ on $\mathbb{R}^d$ is tempered if
there exists $r>0$ such that $\int_{\mathbb{R}^d} (1+|x|^r)^{-1} d\mu(x) < 
+\infty$.)
To underline the fact that the distributions of $W$ are determined by $\Gamma$,
we will write $W_\Gamma$.
The distribution $\Gamma$ is the {\em space correlation} 
of $W_\Gamma$ and $\mu$ is
the {\em spectral measure} of $W_\Gamma$ and~$\Gamma$.

We denote by  $q$, a scalar product  on $S(\mathbb{R}^d)\,$ 
given by the formula: 
$q(\varphi, \psi)= \langle\Gamma,\varphi*\psi_{(s)}\rangle\,,$ where  
$\varphi,\psi \in S(\mathbb{R}^d)$.
In other words, 
such a process $W$ may be called 
{\it associated} with $q$.  

The crucial role in the theory of  stochastic integration with respect
 to  $W_{\Gamma}$  is
played by the Hilbert space $\cH_W\subset S'(\mathbb{R}^d)$ called
the {\em reproducing kernel} Hilbert space of $W_{\Gamma}$.
Namely the space
$\cH_W$ consists of all distributions $\xi\in S'(\mathbb{R}^d)$ for which
there exists a constant $C$
such that $$
|\langle\xi,\psi\rangle|\le C\sqrt{ q(\psi,\psi)},\quad \psi\in S(\mathbb{R}^d).
$$
The norm in $\cH_W$ is given by the formula
$$
|\xi|_{\cH_W}=\sup_{\psi\in S}\frac{|\langle\xi,\psi\rangle|}
{\sqrt{q(\psi,\psi)}}.
$$
\noindent
Let us assume that we require that the stochastic integral should take
values in
a Hilbert space $H\,$ continuously imbedded into $S'(\mathbb{R}^d)\,.$ 
Let $L_{(HS)}(\cH_W,H)$ be
the space of Hilbert-Schmidt operators acting from
$\cH_W$ into $H$. Assume that $\Psi$ is measurable, $(\cF_t)$-adapted,
$L_{HS}(\cH_W,H)$-valued process such that
$$  
E\left( \int_0^t |\Psi(\sigma) |^2_{L_{(HS)}(\cH_W,H)}d\sigma\right)
 <+\infty \quad {\rm for\ all~~} t\ge 0\,.
$$  
Then the stochastic integral
$$
\int_0^t\Psi(\sigma)dW_{\Gamma}(\sigma),\quad t\ge 0
$$
can be defined in a standard way, see \cite{Ito} or \cite{DaPrZa1}. 

As in \cite{PeZa1} and \cite{PeZa2}, we shall use the characterization of
the space $\cH_W$.
\begin{prop} \label{p3.1} (Proposition 1.2, \cite{PeZa1})\\
A distribution $\xi$ belongs to $\cH_W$ if and only if $\xi=\cF(u\mu)$ for
a certain $u\in L^2_{(s)}(\mathbb{R}^d,\mu)$. Moreover, if $\xi=\cF(u\mu)$
and $\eta=\cF(v\mu)$, then 
$$
  \langle\xi,\eta\rangle_{\cH_W} =  
  \langle u,v\rangle_{L^2_{(s)}(\mathbb{R}^d,\mu)}\;,
$$
where $L^2_{(s)}(\mathbb{R}^d,\mu)$ denotes the subspace of 
$L^2_{(s)}(\mathbb{R}^d,\mu;\mathbb{C})$  consisting of all functions $u$
such that $u_{(s)}=u$ and $\langle \cdot ,\cdot\rangle$ denotes products
on particular spaces.
\end{prop} 

\section{Estimates on the interval $[-R,R]$} 

As we have already written, the aim of the paper is to provide conditions under
which the stochastic integral of the form (\ref{eq1}), where  $P=P_\alpha$ is
the fundamental solution of the equation (\ref{eq3}), is a well-defined stochastic
process with values in the space $L^2(\mathbb{R},e^{-|x|}dx)$. 

Let $v$ belong to the space 
$S(\mathbb{R})$ of rapidly decreasing functions on $\mathbb{R}$ 
and be a strictly
positive even function such that $v(x)=e^{-|x|}$ for $|x|\geq 1$. By 
$L^2_v$ we denote the space $L^2(\mathbb{R},vdx)$ which is isomorphic with the
space $L^2(\mathbb{R},e^{-|x|}dx)$. 

\vspace{3mm}
\noindent{\bf Comment:} {\em
 Results obtained in the paper remain true for the function 
$v(x)=(1+|x|^2)^{-\rho}$, with $\rho>\frac{1}{2},~ |x|\geq 1$.}\\[2mm]

Using the notation from section \ref{int}, the stochastic integral 
(\ref{eq1}) reads
\begin{equation}\label{eq9}
 I_\alpha(t) := \int_0^t P_\alpha (t-s) * (b(u(s))\,dW(s)),
\end{equation}
where the convolution means the convolution with respect to the space variable,
that is 
$$ P_\alpha (t-s) * (b(u(s))\,dW(s))(x) =
\int_\mathbb{R} P_\alpha (t-s,y-x)(b(u(s,x))\,dW(s,x))dx. 
$$

We assume that $u$ is an $L^2_v$-valued measurable $({\cal F}_t)$-adapted
process such that 
\begin{equation}\label{eq10}
 \sup_{0\leq t\leq T} E|u(t)|^2_{L^2_v} < +\infty \;.
\end{equation}

Denote by $\mathcal{X}_T$ the space of such processes. In the paper 
$b:~\mathbb{R}\rightarrow\mathbb{R}$ is such a function that the process 
$b(u)$ belongs to the space  $\mathcal{X}_T$. A quite natural example of 
$b$ is any Lipschitz continuous function.\\

We shall use the following\\
{\bf Hypothesis (H):} {\em
~There exists a $\kappa\geq 0$ such that $~\Gamma+\kappa\lambda~$ 
is a non-negative measure, where $\lambda$ denotes Lebesgue measure on
$\mathbb{R}$.}

\vspace{3mm}
\noindent {\bf Comment:}  {\em 
This is known from previous papers (cf\ \cite{PeZa1}, \cite{KaZa2},
 \cite{KaZa3} or \cite{PeZa2}) that the hypothesis (H) is equivalent to
 the assumption:
 there is a constant $\kappa$ such that the measure $\mu + \kappa\delta_0$
 is a positive-definite distribution, where $\mu$ is a spectral measure of the
 noise $W_\Gamma$ and $\delta_0$ is Dirac function. \\
 We recall from section 
 \ref{stoi} that $\Gamma =\cF (\mu)$.
  Next, if $\Gamma$ is a function bounded 
 from below then the  hypothesis (H) holds.
 Additionally, the hypothesis (H) is equivalent to the condition
$$
 \int \limits_{\mathbb{R}} \frac{d\mu(x)}{1+|x|^2} < +\infty\;.
$$
 }

\vspace{1mm}
We define the space $\cH_W^0$ consisting of all distributions of the form 
$\eta=\cF (\psi\mu)$, where $\psi\in S_{(s)}(\mathbb{R})$.

\vspace{3mm}
\noindent {\bf Remarks:} {\em  \vspace{-2mm}
\begin{enumerate}
\item The space $\cH_W^0$  is a dense subspace of $\cH_W$.
\item The following estimate holds
$ \int_\mathbb{R} |\psi(x)|\mu(dx) < +\infty\; ,$ 
where $\psi\in S_{(s)}(\mathbb{R})$.
\item From Prop.~1.3 \cite{PeZa1} and the inclusion  $\cH_W^0\subset H$ 
we obtain $\cH_W^0\subset C_b(\mathbb{R})$.
\end{enumerate}}

\vspace{3mm}
Define $P_\alpha^R(t)$ as follows
$$
 P_\alpha^R(t) := \left\{ \begin{array}{lll}
 P_\alpha(t,|x|) & \mbox{~~for~~} & |x| \leq R\\
 0               & \mbox{~~for~~} & |x| > R
 \end{array} \right.
$$
where $R\in\mathbb{R} $ is finite.

Let us define, analogously like in \cite{PeZa2}, the following operator
\begin{equation} \label{eq11}
 \cK_R(t,u)\,\eta \stackrel{def}{=} P_\alpha^R(t)*(u\eta)
\end{equation}
for $t>0,~u\in L^2_v$ and $\eta\in \cH_W^0$.

Clearly, $$
\cK_R(t,u)\,\eta(x) = \int_\mathbb{R} u(x-y)\eta(x-y\,)P_\alpha^R(t)(dy),
\quad x\in \mathbb{R} \;.$$

In this section we shall show that the operator $\cK_R(t,\cdot)$ 
has an extension to
Hilbert-Schmidt operator from the space $\cH_W$ into $L^2_v$. Additionally,
we will give the estimate
$$  
 |\cK_R(t,u)|_{L_{(HS)}(\cH_W,L^2_v)} \leq C e^R|u|_{L^2_v}\;,
$$  
where $C$ is an appropriate positive constant.
In other words we shall prove that $\cK_R(t,\cdot)$ 
can be uniquelly extended to
linear bounded operator acting from the space $L^2_v$ into the space of
Hilbert-Schmidt operators $L_{(HS)}(\cH_W,L^2_v)$.

First we introduce on the interval $[-R,R]$ the integral analogous
to that defined in~(\ref{eq9}):
\begin{equation} \label{eq12}
 I^R_\alpha(t) :=  \int_0^t P_\alpha^R(t-s)*(b(u(s))dW(s)) =
 \int_0^t \cK_R(t-s,b(u(s)))\,dW(s) \;.
\end{equation}

We can see that the operator $\cK_R$ has to fulfil the condition 
(\ref{eq10}) which guarentees that the integral (\ref{eq12}) is well-defined.
In order to obtain it, the operator $\cK_R$ has to belong to
$L_{(HS)}(\cH_W,H)$. Indeed, if
 $u$ is an $L^2_v$-valued measurable $(\cF_t)$-adapted process fulfilling
condition (\ref{eq10}) then for any $t>0$, the operator-valued process 
$\cK_R(t-s,b(u(s))), s\in(0,t)$, will be adapted and 
will satisfy the required condition
$$ E \left(
\int_0^t |\cK_R(t-s,b(u(s)))|_{L_{(HS)}(\cH_W,H)}\,ds \right)<+\infty\;.
$$ 

Because the process $b(u(t)),~t\in[0,T]$, belongs to the space $\mathcal{X}_T$
for any $u\in\mathcal{X}_T$, in the sequel we shall consider the operators 
$\cK_R(t,u)$ instead of $\cK_R(t,b(u))$.

\vspace{3mm}
Let us emphasize that for any our function 
$v$ there is a constant $C_v$ such that 
\begin{equation} \label{eq13}
 v(x-z) \leq C_v\, e^R\,v(x), \quad
 \mbox{where}\quad x\in \mathbb{R} \quad\mbox{and}\quad z\in [-R,R]\;.
\end{equation}

\begin{lem}\label{lem1}
Assume that $C_v$ fulfills the estimate (\ref{eq13}). Then 
\begin{equation} \label{eq14}
 (P_\alpha^R(t)*v)(x) = \int_\mathbb{R} v(x-y)\,P_\alpha^R(t)\,dy \leq 
  C_v\, e^R\,v(x), \quad\quad x\in[-R,R]\;.
\end{equation}   
Moreover, for all $~t\geq 0$, $~\psi\in S(\mathbb{R}^d)$, the convolution 
$~P_\alpha^R(t)*\psi \in L^2_v~$ and 
\begin{equation} \label{eq15}
 |P_\alpha^R(t)*\psi|_{L^2_v} \leq C_v\,e^R\,|\psi|^2_{L^2_v}\;.
\end{equation}
\end{lem}

\proof{From the property (\ref{eq5}) we obtain
$$(P_\alpha^R(t)*v)(x) = \int_{-R}^{R} v(x-y)\,P_\alpha^R(t)\,dy \leq 
 C_v\, e^R\,v(x) \,\int_{-R}^{R} P_\alpha^R(t)\,dy \leq 
 C_v\, e^R\,v(x) \;,$$
 which proves (\ref{eq14}).
 
 In order to prove (\ref{eq15}) it is enough to write explicitely the 
 norm $|P_\alpha^R(t)*\psi|_{L^2_v}$ for $\psi\in S(\mathbb{R}^d)$
 and then use the estimate (\ref{eq14}).
\hfill $\blacksquare$}\\

\begin{cor}\label{cor2} 
 For any $~t\geq 0~$, there is a unique operator $\cP_\alpha^R(t)\in
 L(L^2_v,L^2_v)$ such that for every $\psi\in S(\mathbb{R})$,
 $$\cP_\alpha^R(t)\,\psi= P_\alpha^R(t) * \psi\, .$$
Additionally, there is a constant $C$ such that 
\begin{equation} \label{eq16}
 |\cP_\alpha^R(t)|_{L(L^2_v,L^2_v)} \leq C\,e^R\;.
\end{equation}
\end{cor}

\vspace{4mm}
From the definition (\ref{eq11}) of the operator $\cK_R$ and the estimation 
(\ref{eq16}), we obtain
$$  
 |\cK_R(t,u)\,\eta|_{L^2_v} \leq C\,e^R\ |u\eta|_{L^2_v}
$$  
for all $~t\geq 0~$, $u\in L^2_v$ and $\eta\in\cH_W^0$.

\vspace{5mm} 
Now, let us assume that the measure $\mu$ satisfies
\begin{equation} \label{eq17}
 \int_\mathbb{R} \frac{d\mu(x)}{1+|x|^2} < +\infty \;.
\end{equation}
We want to show that under the condition  (\ref{eq17}) for any $t>0$,
the operator $\cK_R(t,\cdot)$ has an extension to a bounded linear operator from
the space $L^2_v$ into the space of Hilbert-Schmidt operators
$L_{(HS)}(\cH_W,L^2_v)$.

\vspace{3mm}
\noindent{\bf Remark:} {\em
 Let us note that the hypothesis (H) and the condition (\ref{eq17}) are
  equivalent. This fact comes from Theorem 2, \cite{KaZa2}. }

\vspace{3mm}
\begin{lem} \label{lem2} 
 Assume that $u\in C_b(\mathbb{R})$ and $\{f_k \}\subset \cH_W^0$
 is an orthonormal basis of the space $\cH_W$.
 Then, for any $t>0$, holds
\begin{equation} \label{eq18}
 \sum_{k=1}^{+\infty} |\cK_R(t,u)\,f_k|^2_{L^2_v}  = 
 \int_\mathbb{R}\int_\mathbb{R} 
 |\cF(P_\alpha^R(t)(x-\cdot)u)(y)|^2 \mu(dy)v(x)dx  \;.
\end{equation}
\end{lem}

\vspace{5mm}
\noindent{\bf Comment:} {\em
Lemma \ref{lem2} is formulated for $u\in C_b(\mathbb{R})$, not for $u\in L^2_v$.
This trick provides the sense of the right hand side of (\ref{eq18}).
The main result of this section, Lemma 5, will be formulated for $u\in L^2_v$.
It will be  possible because $C_b(\mathbb{R})$ is a dense subspace of 
the space $L^2_v$.
}\\

\proof{In order to prove (\ref{eq18}) it is enough to rewrite the left hand 
side of (\ref{eq18}), analogously like in \cite{PeZa2}, using the definition
(\ref{eq11}) of the operator $\cK_R$. Next, we have to use properties 
of the Fourier transform and the convolution $P_\alpha^R(t)*uf_k$.
\hfill $\blacksquare$}\\


\begin{lem} \label{lem3} 
 If the spectral measure $\mu$ of the process $W_\Gamma$ fulfills the
condition (\ref{eq17}) then the operators $\cK_R(s,1),~s\geq 0$, acting 
 from the space $\cH_W$
  into $L^2_v$, are Hilbert-Schmidt operators and moreover 
 $$ \int_0^t |\cK_R(s,1)|^2_{L_{HS}(\cH_W,L^2_v)}\,ds <+\infty $$
 for every $t\geq 0$. 
\end{lem}

\vspace{5mm} 
\proof{ Because lemmas analogous to Lemma \ref{lem3} 
 has already been formulated in \cite{KaZa2} and \cite{PeZa2}, 
 we write here only a sketch of the proof. 
 
 We have to obtain the suitable estimate of the right hand side of (\ref{eq18}),
 in the case when the function $u\equiv 1$, in terms of the condition
 (\ref{eq17}). Particularly, we have to estimate the term 
 $|\cF(P_\alpha^R(s)(y))|^2$. In our considerations we shall use formulas 
 (\ref{eq6}),(\ref{eq7}) and Corollary~\ref{cor1}. For symmetry of 
 $P_\alpha^R$ let us consider $y\geq 0$. Then, we may write 
\begin{eqnarray*} 
 \int_{\mathbb{R}} |\cF(P_\alpha^R(s)(y))|^2 \mu(dy) & \leq &
    2\,s^{-\alpha/2} \int_0^{+\infty} 
    |\cF(P_\alpha (ys^{-\alpha/2}))|^2 \mu(dy) \\
 &=&  2\,s^{-\alpha/2} \int_0^{+\infty} 
 \left[\exp \left((ys^{-\alpha/2})^{2/\alpha}
 \cos \frac{\pi}{2}\right)\right]^2 \mu(dy)\;,
\end{eqnarray*}
where $1<\alpha<2$.

For every $1<\alpha<2$ and fixed $s\geq 0$, we may choose a constant
$C_\alpha(s)$ that
$$ \frac{1}{\exp (\frac{a}{s} y^{2/\alpha})^2} \leq
   \frac{C_\alpha(s)}{1+|y|^2}\;,
$$
where $a=|\cos\frac{\pi}{\alpha}|$.

Additionally, we may choose finite constant $C(t)$ that
$$ \int_0^{t} \frac{s^{-\alpha/2}}{\exp (\frac{a}{s} y^{2/\alpha})^2} ds
 \leq \frac{C(t)}{1+|y|^2}\;, \quad y\in\mathbb{R} \;.
$$
Hence, the above estimates give the thesis.
\hfill $\blacksquare$}\\

\vspace{2mm}
\begin{lem} \label{lem4} 
Assume that measure $\mu$ satisfies condition (\ref{eq17}). Then there exists a
constant $C$ such that 
\begin{equation} \label{eq19}
 |\cK_R(t,1)|^2_{L_{HS}(\cH_W,L^2_v)} \leq 
 C_\alpha\left( \int_\mathbb{R} v(x)dx
 \right) \int_\mathbb{R} \frac{\mu (dy)}{1+|y|^2}\;,~~t\geq0\;.
\end{equation}
\end{lem}

\vspace{2mm}
\begin{lem} \label{lem5} 
Assume that the hypothesis (H) holds. Then the operator $\cK_R(t,u)$, for 
$t\geq 0$, $u\in L^2_v$, is a Hilbert-Schmidt operator acting from the space 
$\cH_W$ into $L^2_v$.

\noindent
Additionally, the following estimate is true
\begin{equation} \label{eq20}
 |\cK_R(t,u)|^2_{L_{HS}(\cH_W,L^2_v)} \leq C\, e^R\, |u|^2_{L^2_v} \;,
\end{equation}
where $C$ is an appropriate positive constant.
\end{lem}

\proof{Because proof follows the proof of Lemma 3.3 of \cite{PeZa2} and 
is technical, we present only a sketch of it. We shall use the following
auxiliary result.
\begin{prop} \label{prn}
(Proposition 3.2, \cite{PeZa2})\\
Let $\mu$ be the spectral measure of $W_\Gamma$, 
and let $\delta_0$ be the Dirac
distribution. Condition (H) holds if and only if there exists $\kappa\geq 0$
such that for $N\in\mathbb{N}$, the Fourier transforms of the measures 
$$ \mu_{N,k}(dy) := e^{-|y|^2/N} (\mu(dy)+\kappa\delta_0(dy))\;,
$$ 
are non-negative functions.
\end{prop}
Let us note that now the relationship $\Gamma_{N,k}=\cF(\mu_{N,k})$, 
analogous to that $\Gamma=\cF(\mu)$ is true.

The proof is conducted for functions $u\in C_b(\mathbb{R})$. Because the space 
$C_b$ is dense in $L^2_v$ and the operator $\cK_R(t,u)$ is linear with respect
to $u$, the required estimation obtained for $u\in C_b(\mathbb{R})$ will be
automatically extended to $u\in L^2_v$. 

The proof is done in some steps. First, we consider the case $k=0$ in the above
proposition. We rewrite the formula (\ref{eq18}) from Lemma 2 and express the
integral with respect to the measure $\mu$ like the limit of integrals with
respect to the measures $\mu_{N,0}$ from Proposition \ref{prn}.
Using the formula $\Gamma_{N,0}=\cF(\mu_{N,0})$ and Fubini's theorem, there is
possible to write the series $\sum_{k=1}^{+\infty}|\cK_R(t,u)\,f_k|^2_{L^2_v}$
like integrals expressed in terms of $\Gamma_{N,0}$.  Particularly, it may be
done for $u=1$.

Next, using the estimate (\ref{eq13}), we have the following estimate for all 
 $u\in C_b(\mathbb{R})$ and $t\geq 0$:
$$ \sum_{k=1}^{+\infty}|\cK(t,u)\,f_k|^2_{L^2_v} \leq
 C_v\,e^R\,|u|^2_{L^2_v}\left( \int_{\mathbb{R}} v(x)dx\right)^{-1} \;
 \sum_{k=1}^{+\infty}|\cK(t,1)\,f_k|^2_{L^2_v} \;.
$$ 
From the estimate (\ref{eq19}), we obtain the required result
(\ref{eq20}).\\[-3mm]

The case $k>0$ deals with extending the results obtained in the case $k=0$.
Now, the measure $\nu=\mu+k\delta_0$ satisfies the hypothesis (H). Beacause the
new measure $\nu$ fulfills the condition (\ref{eq17}), we may write the estimate
\begin{equation} \label{eq21}
  |\cK_R(t,u)|^2_{L_{HS}(\cH_V,L^2_v)} \leq C\,e^R\,|u|^2_{L^2_v}
\end{equation}\\[-3mm]
for all $t\geq 0$, $u\in C_b(\mathbb{R})$.
Here $\cH_V$ denotes the reproducing kernel of the spatially homogeneous Wiener
process $V$ with the spectral measure $\nu$.

Basing on the estimate (\ref{eq18}) and the inequality 
$$ \int_{\mathbb{R}} |\cF(P_\alpha^R(t)(x-\cdot))(y)|^2\mu(dy) \leq 
   \int_{\mathbb{R}} |\cF(P_\alpha^R(t)(x-\cdot))(y)|^2\nu(dy)\;,
$$
we have
$$ \sum_{k=1}^{+\infty} |\cK_R(t,u)\,f_k|^2_{L^2_v} \leq
   |\cK_R(t,u)|^2_{L_{HS}(\cH_V,L^2_v)} \;.
$$ 
But from (\ref{eq21}) we obtain $\cK_R(t,u)\in L_{HS}(\cH_W,L^2_v)$
and the estimate (\ref{eq20}), what finishes the proof. 
\hfill $\blacksquare$}

\vspace{3mm}
We have the following consequence of Lemma \ref{lem5}.

\begin{cor}\label{cor3} 
Assume that the hypothesis (H) is satisfied. It is possible to choose a positive
constant $C$ such that if $u\in L^2_v$ and $u(x)\geq 1$ for almost all 
$x\in\mathbb{R}$, then 
$$ \sum_{k=1}^{+\infty} |\cK (t,1)\,f_k|^2_{L^2_v} \leq 
   \sum_{k=1}^{+\infty} |\cK (t,u)\,f_k|^2_{L^2_v} +
   C\,e^R\,|\,u\,|^2_{L^2_v}\,,~~ t\geq 0\;.    $$
In the above estimate the set $\{ f_k\}\subset \cH_W^0$ is an arbitrary 
orthonormal basis of the space~$\cH_W$.  
\end{cor}

\section{Passing to the limit $R\rightarrow +\infty$}

As we have already written, the aim of the paper is to give conditions under
which the integral (\ref{eq9}) is well-defined process with values in the space
$L^2_v$. Till now we have done it on the integral $[-R,R]$ for any finite $R$.
Now, we have to extend our results for $R\rightarrow +\infty$.	

\begin{theor}\label{th4} 
Assume that $b$ is Lipschitz continuous function, $P_\alpha$ is defined by 
(\ref{eq4}) and the hypothesis (H) holds. Then the integral $I_\alpha(t)$
given by (\ref{eq9}) is $L^2_v$-valued.
\end{theor}

\proof{We will show that there exists enough large finite 
$\widetilde{R}>0$ such that for any $R\geq\widetilde{R}$ 
(even $R\rightarrow +\infty$) the following estimate holds

\begin{equation} \label{eq22}
 |\cK_R(t,u)|_{L_{HS}(\cH_W,L^2_v)} \leq 
 |\cK_{\widetilde{R}} (t,u)|_{L_{HS}(\cH_W,L^2_v)} + \cM_{\widetilde{R}}\;,
\end{equation} 
where $\cM_{\widetilde{R}}$ is finite. 

In the 
proof of the estimate (\ref{eq22}) we shall use the definition (\ref{eq11})
of the operator $\cK_R$ and the properties of the function $P_\alpha(t,x)$
defined by (\ref{eq4}). Let us recall that for any $R>0,~t>0,~ x\in \mathbb{R},
~u\in L^2_v~$ and $\eta\in\cH_W^0$
$$
 \cK_R(t,u)\,\eta(x) = \int_\mathbb{R} u(x-y)\,\eta(x-y)\,P_\alpha(t)(y)\,dy\;.
$$

Let us notice that $P_\alpha^R (t)=P_\alpha(t,|x|)$ for $|x|\leq R$, i.e.\
is a fundamental solution to the equation (\ref{eq3}). Additionally,
for $x\geq 0$, $P_\alpha(t,|x|)= P_\alpha(t,x)$, where $P_\alpha(t,x)$  is a 
probability density function (recall (\ref{eq4}),(\ref{eq5}) and (\ref{eq6})).
Because of symmetry of $P_\alpha(t,|x|)$, we will study only the case when 
$x\geq 0$. 

In order to prove the estimate (\ref{eq22}) it is enough to show that for "large"
$\widetilde{R}$ 
$$
 \lim_{R\rightarrow +\infty} \int_{\widetilde{R}}^{R} u(x-y)\,\eta(x-y)
 \,P_\alpha(t,y)\,dy = \cM_{\widetilde{R}}\;,
$$
where $\cM_{\widetilde{R}}< +\infty$, for $x\in\mathbb{R}$ and $t\in(0,+\infty)$.

For any $u\in C_b(\mathbb{R})$ and $\eta\in\cH_W^0\subset C_b(\mathbb{R})$, we
may write 
\begin{equation} \label{eq27}
 \int_{\widetilde{R}}^{R} u(x-y)\,\eta(x-y)\,P_\alpha(t,y)\,dy \leq 
 M_u\, M_\eta  \int_{\widetilde{R}}^{R} P_\alpha(t,y)\,dy \;,
\end{equation} 
where $M_u$ and $M_\eta$ are appropriate constants.

According to Fujita's considerations (see also \cite{Mj} or \cite{Sk}),
the asymptotic behaviour of the density $P_\alpha(1,x)$ is as follows:
\begin{equation} \label{eq28}
 P_\alpha(1,|x|)  \stackrel{|x|\rightarrow +\infty}{\longrightarrow}
 \frac{B_\alpha\, |x|^{(\alpha-1)/(2-\alpha)}}
 {\exp [A_\alpha\,|x|^{2/(2-\alpha)}]}\;,
\end{equation} 
where $A_\alpha$ and $B_\alpha$ are positive constants determined by $\alpha$
with $1\leq \alpha <2$.

Now, let us recall and then use the property (\ref{eq6}), that is 
$P_\alpha(t,x)=P_\alpha(1,x t^{-\alpha/2})\,t^{-\alpha/2}$ for 
$t\in (0,+\infty)$, $x\in\mathbb{R}$.

Then (\ref{eq28}) reads
$$
 P_\alpha(t,x) \;\; \stackrel{|x|\rightarrow +\infty}{\longrightarrow} \;\;
 \frac{B_\alpha\,|xt^{-\alpha/2}|^{(\alpha-1)/(2-\alpha)}t^{-\alpha/2}}
 {\exp [A_\alpha\,|xt^{-\alpha/2}|^{2/(2-\alpha)}]}  \;\;
 \stackrel{|x|\rightarrow +\infty}{\longrightarrow} \;\; 0^+
$$
for $t\in (0,+\infty)$, $x\in\mathbb{R}_+$.

Hence, the right hand side of the estimate (\ref{eq27}) may be writen like
$$
  \cM \int_{\widetilde{R}}^{R} P_\alpha(t,y)\,dy 
  \stackrel{|x|\rightarrow +\infty}{\longrightarrow}
  \cM_{\widetilde{R}}\;,
$$ 
where $\cM_{\widetilde{R}} <+\infty$, what proves (\ref{eq21}).

Since the space $C_b(\mathbb{R})$ is dense in $L^2_v$ and the operator $\cK_R$
is linear with respect to $u$, we obtain the estimate (\ref{eq21}) for all
$t\geq 0$ and $u\in L^2_v$.
\hfill $\blacksquare$}

{\bf Acknowledgements}

The author thanks  Professors Y.~Fujita and S.-O.~Londen for valuable 
remarks concerning integrodifferential equations.\\[-11mm]

\end{document}